\documentclass[a4paper,11pt]{article}
\usepackage{amsmath, amssymb, amsthm}
\usepackage{graphicx} 
\usepackage[all]{xy}
\usepackage[latin1]{inputenc}
\usepackage{enumitem}
\usepackage{color}

\begin{document}

\def\R{\mathbb R} \def\SSS{\mathfrak S} \def\SS{\mathcal S}
\def\XX{\mathcal X} \def\YY{\mathcal Y} \def\VGamma{{\rm V}\Gamma}
\def\Aut{{\rm Aut}} \def\Id{{\rm Id}} \def\I{{\rm I}} 
\def\II{{\rm II}} \def\III{{\rm III}} \def\op{{\rm op}}
\def\ZZ{\mathcal Z} \def\supp{{\rm supp}} \def\Im{{\rm Im}}
\def\sou{{\rm sou}} \def\tar{{\rm tar}}


\title{\bf{A simple solution to the word problem for virtual braid groups}}

\author{\textsc{Paolo Bellingeri, Bruno A. Cisneros de la Cruz, Luis Paris}}

\date{\today}

\maketitle

\begin{abstract}
\noindent
We show a simple and easily implementable solution to the word problem for virtual braid groups.
\end{abstract}


\section{Introduction}

Virtual braid groups were introduced by L. Kauffman in his seminal paper on virtual knots and links \cite{Kauff2}.
They can be defined in several ways, such as in terms of Gauss diagrams \cite{BarDan1, Cisne1}, in terms of braids in thickened surfaces \cite{Cisne1}, and in terms of virtual braid diagrams.
The latter will be our starting point of view.

\bigskip\noindent  
A \emph{virtual braid diagram} on $n$ strands is a $n$-tuple $\beta=(b_1, \dots, b_n)$ of smooth paths in the plane $\R^2$ satisfying the following conditions.
\begin{itemize}[itemsep=2pt,parsep=2pt,topsep=2pt]
\item[(a)]
$b_i(0) = (i,0)$ for all $i \in \{1, \dots, n\}$.
\item[(b)]
There exists a permutation $g \in \SSS_n$ such that $b_i(1) = (g(i),1)$ for all $i \in \{1, \dots, n\}$.
\item[(c)]
$(p_2 \circ b_i)(t)=t$ for all $i \in \{1, \dots, n\}$ and all $t \in [0,1]$, where $p_2: \R^2 \to \R$ denotes the projection on the second coordinate. 
\item[(d)]
The $b_i$'s intersect transversely in a finite number of double points, called the \emph{crossings} of the diagram.
\end{itemize}
Each crossing is endowed with one of the following attributes: positive, negative, virtual.
In the figures they are generally indicated as in Figure 1.1.
Let $VBD_n$ be the set of virtual braid diagrams on $n$ strands, and let $\sim$ be the equivalence relation on $VBD_n$ generated by ambient isotopy and the virtual Reidemeister moves depicted in Figure~1.2.
The concatenation of diagrams induces a group structure on $VBD_n/\sim$.
The latter is called \emph{virtual braid group} on $n$ strands, and is denoted by $VB_n$.

\begin{figure}[ht]
\begin{center}
\begin{tabular}{ccc}
\quad \includegraphics[width=0.8cm]{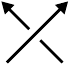} \quad
&
\quad \includegraphics[width=0.8cm]{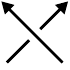} \quad
&
\quad \includegraphics[width=0.8cm]{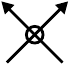} \quad
\\
positive & negative & virtual
\end{tabular}

\bigskip
{\bf Figure 1.1.} Crossings in a virtual braid diagram.
\end{center}
\end{figure}

\begin{figure}[ht]
\begin{center}
\begin{tabular}{cc}
\quad \includegraphics[width=4.8cm]{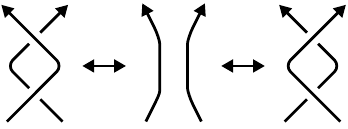} \quad
&
\quad \includegraphics[width=4.4cm]{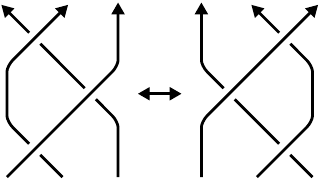} \quad
\end{tabular}

\bigskip
\begin{tabular}{ccc}
\quad \includegraphics[width=2.8cm]{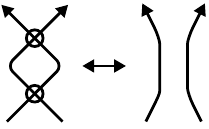} \quad
&
\quad \includegraphics[width=4.4cm]{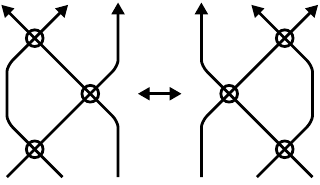} \quad
&
\quad \includegraphics[width=4.4cm]{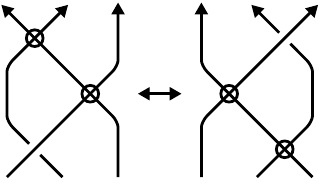} \quad
\end{tabular}

\bigskip
{\bf Figure 1.2.} Virtual Reidemeister moves.
\end{center}
\end{figure}

\bigskip\noindent
It was observed in \cite{Kamad1, Versh1} that $VB_n$ has a presentation with generators $\sigma_1, \dots, \sigma_{n-1}, \tau_1, \dots, \tau_{n-1}$, and relations 
\[ \begin{array}{cl}
\tau_i^2=1 &\ \text{for } 1 \le i \le n-1\\
\sigma_i \sigma_j = \sigma_j \sigma_i\,,\ 
\sigma_i \tau_j = \tau_j \sigma_i\,,  \text{ and }
\tau_i \tau_j = \tau_j \tau_i
&\ \text{for } |i-j| \ge 2\\
\sigma_i \sigma_j \sigma_i = \sigma_j \sigma_i \sigma_j\,,\
\sigma_i \tau_j \tau_i = \tau_j \tau_i \sigma_j\,,\text{ and }
\tau_i \tau_j \tau_i = \tau_j \tau_i \tau_j
& \ \text{for } |i-j|=1
\end{array}\]

\bigskip\noindent
A solution to the word problem for virtual braid groups was shown in \cite{GodPar1}.
However, this solution is quite theoretical and its understanding requires some heavy technical knowledge on Artin groups.
Therefore, it is incomprehensible and useless for most of the potential users, including low dimensional topologists. 
Moreover, its implementation would be difficult.
Our aim here is to show a new solution, which is simpler and easily implementable, and whose understanding does not require any special technical knowledge.
This new solution is in the spirit of the one shown in \cite{GodPar1}, in the sense that one of the main ingredients in its proof is the study of parabolic subgroups in Artin groups. 

\bigskip\noindent
We have not calculated the complexity of this algorithm, as this is probably at least exponential because of the inductive step 3 (see next section). 
Nevertheless, it is quite efficient for a limited number of strands (see the example at the end of Section 2), and, above all, it should be useful to study theoretical questions on $VB_n$ such as the faithfulness of representations of this group in automorphism groups of free groups and/or in linear groups.
Note that the faithfulness of such a representation will immediately provide another, probably faster, solution to the word problem for $VB_n$.

\bigskip\noindent
The Burau representation easily extends to $VB_n$ \cite{Versh1}, but the question whether $VB_n$ is linear or not is still open.
A representation of $VB_n$ in ${\rm Aut}(F_{n+1})$ was independently constructed in \cite{Barda1} and \cite{Mantu1}, but such a representation has recently been proven to be  not faithful for $n \ge 4$ \cite[Proposition 5.3]{Chter1} (see the example at the end of Step 1).
So, we do not know yet any representation on which we can test our algorithm.

\bigskip\noindent
In \cite{Chter1}, Chterental shows a faithful action of $VB_n$ on a set of objects that he calls ``virtual curve diagrams''.
We have some hope to use this action to describe another explicit solution to the word problem for $VB_n$.
But, for now, we do not know any formal definition of this action, and how it could be encoded in an algorithm.

\bigskip\noindent
{\bf Acknowledgments.} 
The research of the first author was partially supported by French grant ANR-11-JS01-002-01.


\section{The algorithm}

Our solution to the word problem for $VB_n$ is divided into four steps.
In Step 1 we define a subgroup $KB_n$ of $VB_n$ and a generating set $\SS$ for $KB_n$, and we show an algorithm (called Algorithm A) which decides whether an element of $VB_n$ belongs to $KB_n$ and, if yes, determines a word over $\SS^{\pm 1}$ which represents this element. 
For $\XX \subset \SS$, we denote by $KB_n(\XX)$ the subgroup of $KB_n$ generated by $\XX$.
The other three steps provide a solution to the word problem for $KB_n(\XX)$ which depends recursively on the cardinality of $\XX$.
Step 2 is the beginning of the induction.
More precisely, the algorithm proposed in Step 2 (called Algorithm B) is a solution to the word problem for $KB_n(\mathcal{X})$ when $\mathcal{X}$ is a full subset of $\mathcal{S}$ (the notion of "full subset" will be also defined in Step 2; for now, the reader just need to know that singletons are full subsets).
In Step 3 we suppose given a solution to the word problem for $KB_n(\XX)$, and, for a given subset $\YY \subset \XX$, we show an algorithm which solves the membership problem for $KB_n(\YY)$ in $KB_n(\XX)$ (Algorithm C).
In Step~4 we show an algorithm which solves the word problem for $KB_n(\XX)$ when $\XX$ is not a full subset, under the assumption that the group $KB_n(\YY)$ has a solvable word problem for any proper subset $\YY$ of $\XX$ (Algorithm~D).
 
\subsection{Step 1}

Recall that $\SSS_n$ denotes the group of permutations of $\{1, \dots, n\}$. 
We denote by $\theta : VB_n \to \SSS_n$ the epimorphism which sends $\sigma_i$ to $1$ and $\tau_i$ to $(i,i+1)$ for all $1 \le i \le n-1$, and by $KB_n$ the kernel of $\theta$.
Note that $\theta$ has a section $\iota : \SSS_n \to VB_n$ which sends $(i,i+1)$ to $\tau_i$ for all $1 \le i \le n-1$, and therefore $VB_n$ is a semi-direct product $VB_n=KB_n \rtimes \SSS_n$. 
The following proposition is proved in Rabenda's master thesis \cite{Raben1} which, unfortunately, is not available anywhere. 
However, its proof can also be found in \cite{BarBel1}.

\bigskip\noindent
{\bf Proposition 2.1} 
(Rabenda \cite{Raben1}).
{\it For $1 \le i < j \le n$ we set 
\[ \begin{array}{rcl}
\delta_{i,j} &=& \tau_i \tau_{i+1} \cdots \tau_{j-2} \sigma_{j-1} \tau_{j-2} \cdots \tau_{i+1} \tau_i\,,\\
\noalign{\smallskip}
\delta_{j,i} &=& \tau_i \tau_{i+1} \cdots \tau_{j-2}\tau_{j-1} \sigma_{j-1} \tau_{j-1} \tau_{j-2} \cdots \tau_{i+1} \tau_i\,.\\
\end{array}\]
Then $KB_n$ has a presentation with generating set
\[
\SS = \{ \delta_{i,j} \mid 1 \le i \neq j \le n\}\,,
\]
and relations 
\[\begin{array}{cl}
\delta_{i,j} \delta_{k,\ell} = \delta_{k,\ell} \delta_{i,j} &\quad \text{for } i,j,k,\ell \text{ distinct}\\
\noalign{\smallskip}
\delta_{i,j} \delta_{j,k} \delta_{i,j} = \delta_{j,k} \delta_{i,j} \delta_{j,k} &\quad \text{for } i,j,k \text{ distinct}
\end{array}\]}

\bigskip\noindent
The virtual braids $\delta_{i,j}$ and $\delta_{j,i}$ are depicted in Figure 2.1.

\begin{figure}[ht]
\begin{center}
\begin{tabular}{cc}
\quad \includegraphics[width=4cm]{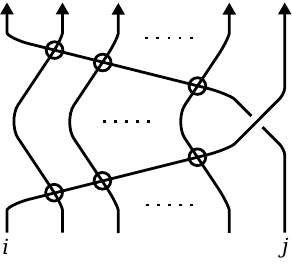} \quad
&
\quad \includegraphics[width=4cm]{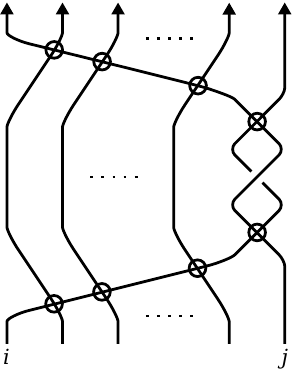} \quad
\\
$\delta_{i,j}$ & $\delta_{j,i}$
\end{tabular}

\bigskip
{\bf Figure 2.1.} Generators for $KB_n$.
\end{center}
\end{figure}

\bigskip\noindent
The following is an important tool in the forthcoming Algorithm A.

\bigskip\noindent
{\bf Lemma 2.2}
(Bardakov, Bellingeri \cite{BarBel1}).
{\it Let $u$ be a word over $\{\tau_1, \dots, \tau_{n-1}\}$, let $\bar u$ be the element of $VB_n$ represented by $u$, and let $i,j \in \{1, \dots, n\}$, $i \neq j$.
Then $\bar u \delta_{i,j} \bar u^{-1} = \delta_{i',j'}$, where $i' = \theta(\bar u)(i)$ and $j'=\theta(\bar u)(j)$.}

\bigskip\noindent
Note that $\tau_i^{-1} = \tau_i$, since $\tau_i^2=1$, for all $i \in \{1, \dots, n-1\}$.
Hence, the letters $\tau_1^{-1}, \dots, \tau_{n-1}^{-1}$ are not needed in the above lemma and below.

\bigskip\noindent
Now, we give an algorithm which, given a word $u$ over $\{ \sigma_1^{\pm 1}, \dots, \sigma_{n-1}^{\pm 1}, \tau_1, \dots, \tau_{n-1}\}$, decides whether the element $\bar u $ of $VB_n$ represented by $u$ belongs to $KB_n$.
If yes, it also determines a word $u'$ over $\SS^{\pm 1}=\{ \delta_{i,j}^{\pm} \mid 1 \le i \neq j \le n\}$ which represents $\bar u$.
The fact that this algorithm is correct follows from Lemma 2.2.

\bigskip\noindent
{\bf Algorithm A.}
Let $u$ be a word over $\{ \sigma_1^{\pm 1}, \dots, \sigma_{n-1}^{\pm 1}, \tau_1, \dots, \tau_{n-1}\}$.
We write $u$ in the form
\[
u= v_0 \sigma_{i_1}^{\varepsilon_1} v_1 \cdots v_{\ell-1} \sigma_{i_\ell}^{\varepsilon_\ell} v_\ell\,,
\]
where $v_0,v_1, \dots, v_\ell$ are words over $\{\tau_1, \dots, \tau_{n-1}\}$, and $\varepsilon_1, \dots, \varepsilon_\ell \in \{ \pm 1\}$.
On the other hand, for a word $v=\tau_{j_1} \cdots \tau_{j_k}$ over $\{\tau_1, \dots, \tau_{n-1}\}$, we set $\theta(v) = (j_1,j_1+1) \cdots (j_k,j_k+1) \in \SSS_n$.
Note that $\theta(\bar u) = \theta(v_0)\, \theta(v_1) \cdots \theta(v_\ell)$.
If $\theta(\bar u) \neq 1$, then $\bar u \not \in KB_n$.
If $\theta(\bar u) = 1$, then $\bar u \in KB_n$, and $\bar u$ is represented by
\[
u' = \delta_{a_1,b_1}^{\varepsilon_1} \delta_{a_2,b_2}^{\varepsilon_2} \cdots \delta_{a_\ell,b_\ell}^{\varepsilon_\ell}\,,
\]
where 
\[
a_k = \theta(v_0 \cdots v_{k-1})(i_k) \text{ and } b_k = \theta(v_0 \cdots v_{k-1})(i_k+1)
\]
for all $k\in\{1, \dots, \ell\}$.

\bigskip\noindent
{\bf Example.} 
In \cite{Chter1} it was proven that the Bardakov-Manturov representation of $VB_n$ in $\Aut(F_{n+1})$ (see for instance \cite{Barda1} for the definition) is not faithful, showing that the element $\omega= (\tau_3 \sigma_2 \tau_1 \sigma_2^{-1})^3$ is non-trivial in  $VB_4$  while  the corresponding automorphism of $F_5$ is trivial.
In \cite{Chter1} the non-triviality of $\omega$ is shown by means of an action on some curve diagrams, but this fact can be easily checked with Algorithm A.
Indeed, $\theta(\omega)=((3,4)(1,2))^3 = (3,4)(1,2) \not= 1$, hence $\omega \not=1$.

\subsection{Step 2}

Let $S$ be a finite set.  
A \emph{Coxeter matrix} over $S$ is a square matrix $M=(m_{s,t})_{s,t \in S}$, indexed by the elements of $S$, such that $m_{s,s}=1$ for all $s \in S$, and $m_{s,t} = m_{t,s} \in \{2,3,4, \dots\} \cup \{ \infty\}$ for all $s,t \in S$, $s \neq t$.
We represent this Coxeter matrix with a labelled graph $\Gamma= \Gamma_M$, called \emph{Coxeter diagram}. 
The set of vertices of $\Gamma$ is $S$.
Two vertices $s,t \in S$ are connected by an edge labelled by $m_{s,t}$ if $m_{s,t} \neq \infty$.

\bigskip\noindent
If $a,b$ are two letters and $m$ is an integer $\ge 2$, we set $\langle a,b \rangle^m = (ab)^{\frac{m}{2}}$ if $m$ is even, and $\langle a,b\rangle^m = (ab)^{\frac{m-1}{2}}a$ if $m$ id odd.
In other words, $\langle a,b \rangle ^m$ denotes the word $aba \cdots$ of length $m$.
The \emph{Artin group} of $\Gamma$ is the group $A=A(\Gamma)$ defined by the following presentation.
\[
A = \langle S \mid \langle s, t \rangle ^{m_{s,t}} = \langle t, s \rangle ^{m_{s,t}} \text{ for all } s,t \in S,\ s \neq t \text{ and } m_{s,t} \neq \infty \rangle\,.
\]
The \emph{Coxeter group} of $\Gamma$, denoted by $W=W(\Gamma)$, is the quotient of $A$ by the relations $s^2=1$, $s \in S$.

\bigskip\noindent
{\bf Example.}
Let $\VGamma_n$ be the Coxeter diagram defined as follows.
The set of vertices of $\VGamma_n$ is $\SS$. 
If $i,j,k,\ell \in \{1, \dots, n\}$ are distinct, then $\delta_{i,j}$ and $\delta_{k,\ell}$ are connected by an edge labelled by $2$. 
If $i,j,k \in \{1, \dots, n\}$ are distinct, then $\delta_{i,j}$ and $\delta_{j,k}$ are connected by an edge labelled by $3$.
There is no other edge in $\VGamma_n$.
Then, by Proposition 2.1, $KB_n$ is isomorphic to $A(\VGamma_n)$.

\bigskip\noindent
Let $\Gamma$ be a Coxeter diagram.
For $X \subset S$, we denote by $\Gamma_X$ the subdiagram of $\Gamma$ spanned by $X$, by $A_X$ the subgroup of $A=A(\Gamma)$ generated by $X$, and by $W_X$ the subgroup of $W=W(\Gamma)$ generated by $X$.
By \cite{Lek1}, $A_X$ is the Artin group of $\Gamma_X$, and, by \cite{Bourb1}, $W_X$ is the Coxeter group of $\Gamma_X$.

\bigskip\noindent
For $\XX \subset \SS$, we denote by $KB_n(\XX)$ the subgroup of $KB_n$ generated by $\XX$. 
By the above, $KB_n(\XX)$ has a presentation with generating set $\XX$ and relations 
\begin{itemize}[itemsep=2pt,parsep=2pt,topsep=2pt]
\item
$st=ts$ if $s$ and $t$ are connected in $\VGamma_n$ by an edge labelled by $2$,
\item
$sts=tst$ if $s$ and $t$ are connected in $\VGamma_n$ by an edge labelled by $3$.
\end{itemize}

\bigskip\noindent
{\bf Definition.}
We say that a subset $\XX$ of $\SS$ is \emph{full} if any two distinct elements $s,t$ of $\XX$ are connected by an edge of $\VGamma_n$.
Recall that the aim of Step 2 is to give a solution to the word problem for $KB_n(\XX)$ when $\XX$ is full.

\bigskip\noindent
We denote by $F_n = F(x_1, \dots, x_n)$ the free group of rank $n$ freely generated by $x_1, \dots, x_n$.
For $i,j \in \{1, \dots, n\}$, $i \neq j$, we define $\varphi_{i,j} \in \Aut(F_n)$ by
\[
\varphi_{i,j} (x_i) = x_i x_j x_i^{-1}\,,\ \varphi_{i,j} (x_j) = x_i\,, \text{ and } \varphi_{i,j} (x_k) = x_k \text{ for } k \not\in \{i,j\}\,.
\]
It is easily checked from the presentation in Proposition 2.1 that the map $\SS \to \Aut(F_n)$, $\delta_{i,j} \mapsto \varphi_{i,j}$, induces a representation $\varphi : KB_n \to \Aut(F_n)$.
For $\XX \subset \SS$, we denote by $\varphi_\XX: KB_n(\XX) \to \Aut(F_n)$ the restriction of $\varphi$ to $KB_n(\XX)$.
The following will be proved in Section 3.

\bigskip\noindent
{\bf Proposition 2.3.}
{\it If $\XX$ is a full subset of $\SS$, then $\varphi_\XX : KB_n(\XX) \to \Aut(F_n)$ is faithful.}

\bigskip\noindent
{\bf Notation.}
From now on, if $u$ is a word over $\SS^{\pm 1}$, then $\bar u$ will denote the element of $KB_n$ represented by $u$.

\bigskip\noindent
{\bf Algorithm B.}
Let $\XX$ be a full subset of $\SS$, and let $u=s_1^{\varepsilon_1} \cdots s_\ell^{\varepsilon_\ell}$ be a word over $\XX^{\pm 1}$.
We have $\varphi_\XX(\bar u) = \varphi_\XX(s_1)^{\varepsilon_1} \cdots \varphi_\XX(s_\ell)^{\varepsilon_\ell}$.
If $\varphi(\bar u) = \Id$, then $\bar u=1$.
Otherwise, $\bar u \neq 1$.

\subsection{Step 3}

Let $G$ be a group, and let $H$ be a subgroup of $G$.
A solution to the \emph{membership problem} for $H$ in $G$ is an algorithm which, given $g \in G$, decides whether $g$ belongs to $H$ or not.
In the present step we will assume that $KB_n(\XX)$ has a solution to the word problem, and, from this solution, we will give a solution to the membership problem for $KB_n(\YY)$ in $KB_n(\XX)$, for $\YY \subset \XX$. 
Furthermore, if the tested element belongs to $KB_n(\YY)$, then this algorithm will determine a word over $\YY^{\pm 1}$ which represents this element. 

\bigskip\noindent
Let $u$ be a word over $\SS$.
(Remark: here the alphabet is $\SS$, and not $\SS^{\pm 1}$.)
\begin{itemize}[itemsep=2pt,parsep=2pt,topsep=2pt]
\item
Suppose that $u$ is written in the form $u_1 ss u_2$, where $u_1,u_2$ are words over $\SS$ and $s$ is an element of $\SS$. 
Then we say that $u' = u_1u_2$ is obtained from $u$ by an \emph{$M$-operation of type $\I$}.
\item
Suppose that $u$ is written in the form $u_1 st u_2$, where $u_1,u_2$ are words over $\SS$ and $s,t$ are two elements of $\SS$ connected by an edge labelled by $2$.
Then we say that $u' = u_1 ts u_2$ is obtained from $u$ by an \emph{$M$-operation of type $\II^{(2)}$}.
\item
Suppose that $u$ is written in the form $u_1 sts u_2$, where $u_1,u_2$ are words over $\SS$ and $s,t$ are two elements of $\SS$ connected by an edge labelled by $3$.
Then we say that $u' = u_1 tst u_2$ is obtained from $u$ by an \emph{$M$-operation of type $\II^{(3)}$}.
\end{itemize}
Let $\YY$ be a subset of $\SS$.
\begin{itemize}[itemsep=2pt,parsep=2pt,topsep=2pt]
\item
Suppose that $u$ is written in the form $tu'$, where $u'$ is a word over $\SS$ and $t$ is an element of $\YY$.
Then we say that $u'$ is obtained from $u$ by an \emph{$M$-operation of type $\III_\YY$}.
\end{itemize}

\bigskip\noindent
We say that $u$ is \emph{$M$-reduced} (resp. \emph{$M_\YY$-reduced}) if its length cannot be shortened by $M$-operations of type $\I,\II^{(2)}, \II^{(3)}$ (resp. of type $\I,\II^{(2)}, \II^{(3)}, \III_\YY$).
An \emph{$M$-reduction} (resp. \emph{$M_\YY$-reduction}) of $u$ is an $M$-reduced word (resp. $M_\YY$-reduced word) obtained from $u$ by $M$-operations (resp. $M_\YY$-operations).
We can easily enumerate all the words obtained from $u$ by $M$-operations (resp. $M_\YY$-operations), hence we can effectively determine an $M$-reduction and/or an $M_\YY$-reduction of $u$. 

\bigskip\noindent
Let $\YY$ be a subset of $\SS$.
From a word $u =s_1^{\varepsilon_1} \cdots s_\ell^{\varepsilon_\ell}$ over $\SS^{\pm 1}$, we construct a word $\pi_\YY(u)$ over $\YY^{\pm 1}$ as follows.
\begin{itemize}[itemsep=2pt,parsep=2pt,topsep=2pt]
\item
For $i \in \{0,1, \dots, \ell\}$ we set $u_i^+=s_1 \cdots s_i$ (as ever, $u_0^+$ is the identity).
\item
For $i \in \{0,1, \dots, \ell\}$ we calculate an $M_\YY$-reduction $v_i^+$ of $u_i^+$.
\item
For a word $v=t_1 \cdots t_k$ over $\SS$, we denote by $\op(v) = t_k \cdots t_1$ the \emph{anacycle} of $v$.
Let $i \in \{1, \dots, \ell\}$.
If $\varepsilon_i=1$, we set $w_i^+ = v_{i-1}^+\cdot s_i\cdot \op(v_{i-1}^+)$.
If $\varepsilon_i=-1$, we set $w_i^+ = v_{i}^+\cdot s_i \cdot\op(v_{i}^+)$.
\item
For all $i \in \{1, \dots, \ell\}$ we calculate an $M$-reduction $r_i$ of $w_i^+$.
\item
If $r_i$ is of length $1$ and $r_i \in \YY$, we set $T_i=r_i^{\varepsilon_i}$.
Otherwise we set $T_i=1$.
\item
We set $\pi_\YY(u) = T_1 T_2 \cdots T_\ell$.
\end{itemize}

\bigskip\noindent
The proof of the following is given in Section 4. 

\bigskip\noindent
{\bf Proposition 2.4.}
{\it Let $\YY$ be a subset of $\SS$.
Let $u,v$ be two words over $\SS^{\pm 1}$.
If $\bar u = \bar v$, then $\overline{\pi_\YY(u)} = \overline{\pi_\YY(v)}$.
Moreover, we have $\bar u \in KB_n(\YY)$ if and only if $\bar u = \overline{\pi_\YY(u)}$.}

\bigskip\noindent
{\bf Algorithm C.}
Take two subsets $\XX$ and $\YY$ of $\SS$ such that $\YY \subset \XX$, and assume given a solution to the word problem for $KB_n(\XX)$.
Let $u$ be a word over $\XX^{\pm 1}$.
We calculate $v=\pi_\YY(u)$.
If $\overline{uv^{-1}} \neq 1$, then $\bar u \not\in KB_n(\YY)$.
If $\overline{uv^{-1}} = 1$, then $\bar u \in KB_n(\YY)$ and $v$ is a word over $\YY^{\pm 1}$ which represents the same element as $u$.

\bigskip\noindent
We can use Algorithm C to show that the representation $\varphi: KB_n \to \Aut(F_n)$ of Step 2 is not faithful.
Indeed, let $\alpha = \delta_{1,3} \delta_{3,2} \delta_{3,1}$ and $\beta = \delta_{2,3} \delta_{1,3} \delta_{3,2}$.
A direct calculation shows that $\varphi (\alpha) = \varphi (\beta)$.
Now, set $\XX = \SS$ and $\YY = \{ \delta_{1,3}, \delta_{3,2}, \delta_{3,1} \}$.
We have $\pi_\YY (\delta_{1,3} \delta_{3,2} \delta_{3,1}) = \delta_{1,3} \delta_{3,2} \delta_{3,1}$, hence $\alpha \in KB_n(\YY)$, and we have $\pi_\YY(\delta_{2,3} \delta_{1,3} \delta_{3,2})=1$ and $\beta \neq 1$, hence $\beta \not\in KB_n(\YY)$.
So, $\alpha \neq \beta$.


\subsection{Step 4}

Now, we assume that $\XX$ is a non-full subset of $\SS$, and that we have a solution to the word problem for $KB_n(\YY)$ for any proper subset $\YY$ of $\XX$ (induction hypothesis). 
We can and do choose two proper subsets $\XX_1, \XX_2 \subset \XX$ satisfying the following properties.
\begin{itemize}[itemsep=2pt,parsep=2pt,topsep=2pt]
\item[(a)]
$\XX=\XX_1 \cup \XX_2$.
\item[(b)]
Let $\XX_0 = \XX_1 \cap \XX_2$.
There is no edge in $\VGamma_n$ connecting an element of $\XX_1 \setminus \XX_0$ to an element of $\XX_2 \setminus \XX_0$.
\end{itemize}
It is easily seen from the presentations of the $KB_n(\XX_i)$'s given in Step~2 that we have the amalgamated product 
\[
KB_n(\XX) = KB_n(\XX_1) *_{KB_n(\XX_0)} KB_n(\XX_2)\,.
\]

\bigskip\noindent
Our last algorithm is based on the following result.
This is well-known and can be found for instance in \cite[Chap. 5.2]{Serre1}.

\bigskip\noindent
{\bf Proposition 2.5.}
{\it Let $A_1 *_B A_2$ be an amalgamated product of groups. 
Let $g_1, \dots, g_\ell$ be a sequence of elements of $A_1 \sqcup A_2$ different from $1$ and satisfying the following condition:  \begin{itemize}
\item[]
if $g_i \in A_1$ (resp. $g_i \in A_2$), then $g_{i+1} \in A_2 \setminus B$ (resp.  $g_{i+1} \in A_1 \setminus B$), for all $i \in \{1, \dots, \ell-1\}$.
\end{itemize}
Then $g_1 g_2 \cdots g_\ell$ is different from $1$ in $A_1 *_B A_2$.}

\bigskip\noindent
{\bf Algorithm D.}
Let $u$ be a word over $\XX^{\pm 1}$.
We write $u$ in the form $u_1u_2 \cdots u_\ell$, where
\begin{itemize}[itemsep=2pt,parsep=2pt,topsep=2pt]
\item
$u_i$ is either a word over $\XX_1^{\pm 1}$, or a word over $\XX_2^{\pm 1}$,
\item
if $u_i$ is a word over $\XX_1^{\pm 1}$ (resp. over $\XX_2^{\pm 1}$), then $u_{i+1}$ is a word over $\XX_2^{\pm 1}$ (resp. over $\XX_1^{\pm 1}$).
\end{itemize}
We decide whether $\bar u$ is trivial by induction on $\ell$.
Suppose that $\ell=1$ and $u=u_1 \in KB_n(\XX_j)$ ($j \in \{1,2\}$).
Then we apply the solution to the word problem for $KB_n(\XX_j)$ to decide whether $\bar u$ is trivial or not.
Suppose that $\ell \ge 2$.
For all $i$ we set $v_i=\pi_{\XX_0}(u_i)$.
If $\overline{u_iv_i^{-1}} \neq 1$ for all $i$, then $\bar u \neq 1$.
Suppose that there exists $i \in \{1, \dots, \ell\}$ such that $\overline{u_iv_i^{-1}} = 1$.
Let $u_i' = v_1u_2$ if $i=1$, $u_i' = u_{\ell-1} v_\ell$ if $i=\ell$, and $u_{i}'=u_{i-1} v_i u_{i+1}$ if $2 \le i \le \ell-1$.
Set $v=u_1 \cdots u_{i-2} u_{i}' u_{i+2} \cdots u_\ell$.
Then $\bar u = \bar v$ and, by induction, we can decide whether $v$ represents $1$ or not. 

\subsection{Example}

In order to illustrate our solution to the word problem for $KB_n$, we turn now to give a more detailed and efficient version of the algorithm for the group $KB_4$.
We start with the following observation.

\bigskip\noindent
{\bf Remark.}
For $\XX \subset \SS$, we denote by $\VGamma_n(\XX)$ the full subgraph of $\VGamma_n$ spanned by $\XX$.
Let $\XX, \YY$ be two subsets of $\SS$.
Note that an injective morphism of Coxeter graphs $\VGamma_n(\YY) \hookrightarrow \VGamma_n (\XX)$ induces an injective homomorphism $KB_n (\YY) \hookrightarrow KB_n(\XX)$.
So, if we have a solution to the word problem for $KB_n(\XX)$, then such a morphism would determine a solution to the word problem for $KB_n(\YY)$.

\bigskip\noindent
The Coxeter graph $\VGamma_4$ is depicted in Figure 2.2.
Our convention in this figure is that a full edge is labelled by $3$ and a dotted edge is labelled by $2$. 
Note that there are two edges that go through ``infinity'', one connecting $\delta_{2,1}$ to $\delta_{4,3}$, and one connecting $\delta_{1,4}$ to $\delta_{3,2}$.

\begin{figure}[ht]
\begin{center}
\includegraphics[width=4.8cm]{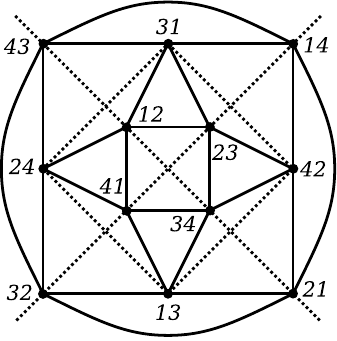}

\bigskip
{\bf Figure 2.2.} Coxeter graph $\VGamma_4$.
\end{center}
\end{figure}

\bigskip\noindent
Consider the following subsets of $\SS$.
\begin{itemize}[itemsep=2pt,parsep=2pt,topsep=2pt]
\item[]
$\XX(1) = \{ \delta_{1,2}, \delta_{2,3}, \delta_{3,4}, \delta_{4,1}, \delta_{3,1} \}$, $\XX_1(1) = \{ \delta_{1,2}, \delta_{2,3}, \delta_{3,4}, \delta_{4,1} \}$, $\XX_2(1) = \{ \delta_{1,2}, \delta_{2,3}, \delta_{3,1} \}$.
\item[]
$\XX(2) = \XX(1) \cup \{ \delta_{4,2} \}$, $\XX_1(2) = \XX(1)$, $\XX_2(2) = \{ \delta_{4,2}, \delta_{3,4}, \delta_{2,3}, \delta_{3,1} \}$.
\item[]
$\XX(3) = \XX(2) \cup \{ \delta_{1,3} \}$, $\XX_1(3) = \XX(2)$, $\XX_2(3) = \{ \delta_{1,3}, \delta_{4,1}, \delta_{3,4}, \delta_{4,2} \}$.
\item[]
$\XX(4) = \XX(3) \cup \{ \delta_{2,4} \}$, $\XX_1(4) = \XX(3)$, $\XX_2(4) = \{ \delta_{2,4}, \delta_{1,3}, \delta_{ 4,1}, \delta_{1,2}, \delta_{3,1} \}$.
\item[]
$\XX(5) = \XX(4) \cup \{ \delta_{1,4} \}$, $\XX_1(5) = \XX(4)$, $\XX_2(5) = \{ \delta_{1,4}, \delta_{4,2}, \delta_{2,3}, \delta_{3,1} \}$.
\item[]
$\XX(6) = \XX(5) \cup \{ \delta_{2,1} \}$, $\XX_1(6) = \XX(5)$, $\XX_2(6) = \{ \delta_{2,1}, \delta_{1,3}, \delta_{3,4}, \delta_{4,2}, \delta_{1,4} \}$.
\item[]
$\XX(7) = \XX(6) \cup \{ \delta_{3,2} \}$, $\XX_1(7) = \XX(6)$, $\XX_2(7) = \{ \delta_{3,2}, \delta_{2,4}, \delta_{4,1}, \delta_{1,3}, \delta_{2,1}, \delta_{1,4} \}$.
\item[]
$\XX(8) = \XX(7) \cup \{ \delta_{4,3} \} = \SS$, $\XX_1(8) = \XX(7)$, $\XX_2(8) = \{ \delta_{4,3}, \delta_{3,2}, \delta_{2,4}, \delta_{1,2}, \delta_{3,1}, \delta_{1,4}, \delta_{2,1} \}$.
\end{itemize}

\bigskip\noindent
Let $k \in \{1, \dots, 8\}$.
Note that $\XX(k) = \XX_1(k) \cup \XX_2(k)$.
The Coxeter graph $\VGamma_4(\XX(k))$ is depicted in Figure 2.3.
In this figure the elements of $\XX_1(k)$ are represented by punctures, while the elements of $\XX_2(k)$ are represented by small circles.

\begin{figure}[ht]
\begin{center}
\begin{tabular}{cccc}
\parbox[c]{1.2cm}{\includegraphics[width=1cm]{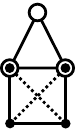}}
&
\parbox[c]{2cm}{\includegraphics[width=1.8cm]{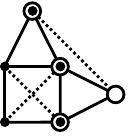}}
&
\parbox[c]{2cm}{\includegraphics[width=1.8cm]{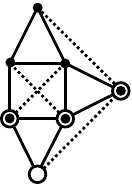}}
&
\parbox[c]{2.8cm}{\includegraphics[width=2.6cm]{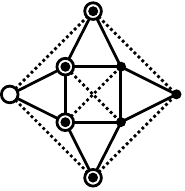}}
\\
$\VGamma_4(\XX(1))$ & $\VGamma_4(\XX(2))$ & $\VGamma_4(\XX(3))$ & $\VGamma_4(\XX(4))$
\end{tabular}

\bigskip
\begin{tabular}{cccc}
\parbox[c]{2.8cm}{\includegraphics[width=2.6cm]{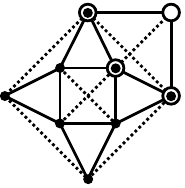}}
&
\parbox[c]{3cm}{\includegraphics[width=2.8cm]{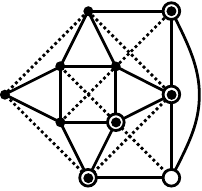}}
&
\parbox[c]{3.4cm}{\includegraphics[width=3.2cm]{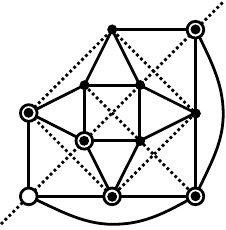}}
&
\parbox[c]{3.4cm}{\includegraphics[width=3.2cm]{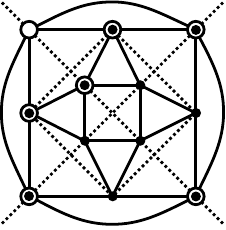}}
\\
$\VGamma_4(\XX(5))$ & $\VGamma_4(\XX(6))$ & $\VGamma_4(\XX(7))$ & $\VGamma_4(\XX(8))$
\end{tabular}

\bigskip
{\bf Figure 2.3.} Coxeter graph $\VGamma_4(\XX(k))$.
\end{center}
\end{figure}

\bigskip\noindent
We solve the word problem for $KB_4 (\XX (k))$ successively for $k=1,2, \dots, 8$, thanks to the following observations. 
Since $\XX(8) =\SS$, this will provide a solution to the word problem for $KB_4$.
\begin{itemize}[itemsep=2pt,parsep=2pt,topsep=2pt]
\item[(1)]
Let $k \in \{1, \dots, 8\}$.
Set $\XX_0(k) = \XX_1(k) \cap \XX_2(k)$.
Observe that there is no edge in $\VGamma_4$ connecting an element of $\XX_1(k) \setminus \XX_0 (k)$ to an element of $\XX_2(k) \setminus \XX_0 (k)$.
Hence, we can solve with Algorithm D the word problem for $KB_4(\XX(k))$ from solutions to the word problem for $KB_4 (\XX_1 (k))$ and for $KB_4 (\XX_2 (k))$.
\item[(2)]
The subsets $\XX_1(1)$ and $\XX_2(1)$ are full, hence we can solve the word problem for $KB_4(\XX_1 (1))$ and for $KB_4 (\XX_2 (1))$ with Algorithm B.
\item[(3)]
Let $k \ge 2$.
On the one hand, we have $\XX_1(k) = \XX(k-1)$.
On the other hand, it is easily seen that there is an injective morphism $\VGamma_4 (\XX_2 (k)) \hookrightarrow \VGamma_4 (\XX (k-1))$.
Hence, by the remark given at the beginning of the subsection, we can solve the word problem for $KB_4 (\XX_1 (k))$ and for $KB_4 (\XX_2 (k))$ from a solution to the word problem for $KB_4 (\XX (k-1))$.
\end{itemize}

\section{Proof of Proposition 2.3}

Recall that $F_n=F(x_1, \dots, x_n)$ denotes the free group of rank $n$ freely generated by $x_1, \dots, x_n$, and that we have a representation $\varphi : KB_n \to \Aut(F_n)$ which sends $\delta_{i,j}$ to $\varphi_{i,j}$, where
\[
\varphi_{i,j} (x_i) = x_i x_j x_i^{-1}\,,\ \varphi_{i,j} (x_j) = x_i\,, \text{ and } \varphi_{i,j} (x_k) = x_k \text{ for } k \not\in \{i,j\}\,.
\]
For $\XX \subset \SS$, we denote by $\varphi_\XX: KB_n(\XX) \to \Aut(F_n)$ the restriction of $\varphi$ to $KB_n(\XX)$.
In this section we prove that $\varphi_\XX$ is faithful if $\XX$ is a full subset of $\SS$.

\bigskip\noindent
Consider the following groups.
\begin{gather*}
B_n=\left\langle \sigma_1, \dots, \sigma_{n-1} \ \left|
\begin{array}{cl}
\sigma_i \sigma_j \sigma_i = \sigma_j \sigma_i \sigma_j &\text{if } |i-j|=1\\ 
\sigma_i \sigma_j = \sigma_j \sigma_i &\text{if } |i-j|\ge 2 
\end{array}\right. \right\rangle\,,
\\
\tilde B_n = \left\langle \sigma_1, \dots, \sigma_n\ \left|
\begin{array}{cl}
\sigma_i \sigma_j \sigma_i = \sigma_j \sigma_i \sigma_j &\text{if } i \equiv j \pm 1 \mod n\\
\sigma_i \sigma_j = \sigma_j \sigma_i &\text{if } i \neq j \text{ and } i \not\equiv j \pm 1 \mod n
\end{array} \right. \right\rangle\,,\quad n \ge 3\,.
\end{gather*}
The group $B_n$ is the classical \emph{braid group}, and $\tilde B_n$ is the \emph{affine braid group}. 

\bigskip\noindent
We define representations $\psi_n : B_n \to \Aut(F_n)$  and $\tilde \psi_n : \tilde B_n  \to \Aut(F_n)$ in the same way as $\varphi$ as follows. 
\begin{gather*}
\psi_n(\sigma_i)(x_i) = x_i x_{i+1} x_i^{-1}\,,\ \psi_n(\sigma_i)(x_{i+1}) = x_i\,,\ \psi_n(\sigma_i)(x_k) = x_k\ \text{if } k \not\in\{i,i+1\}\\
\tilde\psi_n(\sigma_i)(x_i) = x_i x_{i+1} x_i^{-1}\,,\ \tilde\psi_n(\sigma_i)(x_{i+1}) = x_i\,,\ \tilde\psi_n(\sigma_i)(x_k) = x_k\ \text{if } k \not\in\{i,i+1\}\,,\ \text{for }i<n\\
\tilde\psi_n(\sigma_n)(x_n) = x_nx_1x_n^{-1}\,,\ \tilde\psi_n(\sigma_n)(x_1)=x_n\,,\ \tilde\psi_n(\sigma_n)(x_k)=x_k \text{ if } k \not\in \{1,n\}
\end{gather*}
The key of the proof of Proposition 2.3 is the following. 

\bigskip\noindent
{\bf Theorem 3.1}
(Artin \cite{Artin2}, Bellingeri, Bodin \cite{BelBod1}).
{\it The representations $\psi_n: B_n \to \Aut(F_n)$ and $\tilde \psi_n: \tilde B_n \to \Aut(F_n)$ are faithful.}

\bigskip\noindent
The \emph{support} of a generator $\delta_{i,j}$ is defined to be $\supp(\delta_{i,j})=\{i,j\}$.
The \emph{support} of a subset $\XX$ of $\SS$ is $\supp(\XX)=\cup_{s \in \XX} \supp(s)$.
We say that two subsets $\XX_1$ and $\XX_2$ of $\SS$ are \emph{perpendicular}\footnote{This terminology is derived from the theory of Coxeter groups.}
if $\supp(\XX_1) \cap \supp(\XX_2) = \emptyset$.
Note that this condition implies that $\XX_1 \cap \XX_2 = \emptyset$.
More generally, we say that a family $\XX_1, \dots, \XX_\ell$ of subsets of $\SS$ is \emph{perpendicular} if $\supp(\XX_i) \cap \supp(\XX_j) = \emptyset$ for all $i \neq j$.
In that case we write $\XX_1 \cup \cdots \cup \XX_\ell = \XX_1 \boxplus \cdots \boxplus \XX_\ell$.
We say that a subset $\XX$ of $\SS$ is \emph{indecomposable} if it is not the union of two perpendicular nonempty subsets. 
The following observations will be of importance in what follows.

\bigskip\noindent
{\bf Remark.}
Let $\XX_1$ and $\XX_2$ be two perpendicular subsets of $\SS$, and let $\XX = \XX_1 \boxplus \XX_2$.
\begin{itemize}[itemsep=2pt,parsep=2pt,topsep=2pt]
\item[(1)]
$\XX$ is a full subset if and only if $\XX_1$ and $\XX_2$ are both full subsets.
\item[(2)]
$KB_n(\XX) = KB_n(\XX_1) \times KB_n(\XX_2)$.
\end{itemize}
Indeed, if $\delta_{i,j} \in \XX_1$ and $\delta_{k,\ell} \in \XX_2$, then $i,j,k,\ell$ are distinct, and therefore $\delta_{i,j}$ and $\delta_{k,\ell}$ are connected by an edge labelled by $2$, and $\delta_{i,j} \delta_{k,\ell} = \delta_{k,\ell} \delta_{i,j}$.

\bigskip\noindent
{\bf Lemma 3.2.}
{\it Let $\XX_1$ and $\XX_2$ be two perpendicular subsets of $\SS$, and let $\XX = \XX_1 \boxplus \XX_2$.
Then $\varphi_\XX : KB_n(\XX) \to \Aut(F_n)$ is faithful if and only if $\varphi_{\XX_1} : KB_n(\XX_1) \to \Aut(F_n)$ and $\varphi_{\XX_2} : KB_n(\XX_2) \to \Aut(F_n)$ are both faithful.}

\bigskip\noindent
{\bf Proof.}
For $X \subset \{x_1, \dots, x_n\}$, we denote by $F(X)$ the subgroup of $F_n$ generated by $X$.
There is a natural embedding $\iota_X : \Aut(F(X)) \hookrightarrow \Aut(F_n)$ defined by
\[
\iota_X(\alpha)(x_i) = \left\{ \begin{array}{ll}
\alpha(x_i) &\text{if } x_i \in X\\
x_i &\text{otherwise}
\end{array} \right.
\]
Moreover, if $X_1$ and $X_2$ are disjoint subsets of $\{x_1, \dots, x_n\}$, then the homomorphism
\[
\begin{array}{rccc}
(\iota_{X_1} \times \iota_{X_2}) : & \Aut(F(X_1)) \times \Aut(F(X_2)) & \to & \Aut (F_n)\\
& (\alpha_1, \alpha_2) & \mapsto & \iota_{X_1} (\alpha_1)\,\iota_{X_2}(\alpha_2)
\end{array}
\]
is well-defined and injective.
From now on, we will assume $\Aut(F(X))$ to be embedded in $\Aut(F_n)$ via $\iota_X$, for all $X \subset \{x_1, \dots, x_n\}$.

\bigskip\noindent
By abuse of notation, for $\XX \subset \SS$, we will also denote by $\supp(\XX)$ the set $\{x_i \mid i \in \supp(\XX)\}$.
Set $X_1 = \supp(\XX_1)$ and $X_2 = \supp(\XX_2)$.
We have $\Im(\varphi_{\XX_i}) \subset \Aut(F(X_i))$ for $i=1,2$, $X_1 \cap X_2 = \emptyset$, and $KB_n(\XX) = KB_n(\XX_1) \times KB_n(\XX_2)$.
Hence, Lemma 3.2 follows from the following claim whose proof is left to the reader.

\bigskip\noindent
{\it Claim.}
Let $f_1 : G_1 \to H_1$ and $f_2 : G_2 \to H_2$ be two group homomorphisms. 
Let $(f_1 \times f_2) : (G_1 \times G_2) \to (H_1 \times H_2)$ be the homomorphism defined by $(f_1 \times f_2) (u_1,u_2) = (f_1(u_1), f_2(u_2))$.
Then $(f_1 \times f_2)$ is injective if and only if $f_1$ and $f_2$ are both injective.
\qed

\bigskip\noindent
For $2\le m\le n$ we set 
\[
\ZZ_m=\{ \delta_{1,2}, \dots, \delta_{m-1,m}\}\,,\ \tilde \ZZ_m=\{ \delta_{1,2}, \dots, \delta_{m-1,m}, \delta_{m,1}\}\,.
\]
Note that the map $\{\sigma_1, \dots, \sigma_{m-1} \} \to \ZZ_m$, $\sigma_i \mapsto \delta_{i,i+1}$, induces an isomorphism $f_m: B_m \to KB_n(\ZZ_m)$.
This follows from the presentation of $KB_n(\ZZ_m)$ given in Step 2 of Section 2.
Similarly, for $m \ge 3$, the map $\{\sigma_1, \dots, \sigma_m\} \to \tilde \ZZ_m$, $\sigma_i \mapsto \delta_{i,i+1}$ for $1 \le i \le m-1$, $\sigma_m \mapsto \delta_{m,1}$, induces an isomorphism $\tilde f_m: \tilde B_m \to KB_n(\tilde \ZZ_m)$.

\bigskip\noindent
Recall that the symmetric group $\SSS_n$ acts on $\SS$ by $g\,\delta_{i,j} = \delta_{g(i), g(j)}$, and that this action induces an action of $\SSS_n$ on $KB_n$.
On the other hand, there is a natural embedding $\SSS_n \hookrightarrow \Aut(F_n)$, where $g \in \SSS_n$ sends $x_i$ to $x_{g(i)}$ for all $i \in \{1, \dots, n\}$, and this embedding induces by conjugation an action of $\SSS_n$ on $\Aut(F_n)$.
It is easily seen that the homomorphism $\varphi : KB_n \to \Aut(F_n)$ is equivariant under these actions of $\SSS_n$.

\bigskip\noindent
{\bf Lemma 3.3.}
{\it If $\XX$ is a full and indecomposable nonempty subset of $\SS$, then there exist $g \in \SSS_n$ and $m \in \{2, \dots, n\}$ such that either $\XX = g\,\ZZ_m$, or $\XX=g\,\tilde\ZZ_m$ and $m \ge 3$.}

\bigskip\noindent
{\bf Proof.}
An \emph{oriented graph} $\Upsilon$ is the data of two sets, $V(\Upsilon)$, called \emph{set of vertices}, and $E(\Upsilon)$, called \emph{set of arrows}, together with two maps $\sou, \tar : E(\Upsilon) \to V(\Upsilon)$.
We associate an oriented graph $\Upsilon_\XX$ to any subset $\XX$ of $\SS$ as follows.
The set of vertices is $V(\Upsilon_\XX) = \supp (\XX)$, the set of arrows is $E(\Upsilon_\XX) = \XX$, and, for $\delta_{i,j} \in \XX$, we set $\sou (\delta_{i,j})=i$ and $\tar (\delta_{i,j})=j$.
Assume that $\XX$ is a full and indecomposable nonempty subset of $\SS$.
Since $\XX$ is indecomposable, $\Upsilon_\XX$ must be connected.
Since $\XX$ is full, if $s,t \in \XX$ are two different arrows of $\Upsilon_\XX$ with a common vertex, then there exist $i,j,k \in \{1, \dots, n\}$ distinct such that either $s=\delta_{j,i}$ and $t=\delta_{i,k}$, or $s=\delta_{i,j}$ and $t = \delta_{k,i}$.
This implies that $\Upsilon_\XX$ is either an oriented segment, or an oriented cycle with at least $3$ vertices (see Figure 3.1). 
If $\Upsilon_\XX$ is an oriented segment, then there exist $g \in \SSS_n$ and $m \in \{2, \dots, n\}$ such that $\XX=g\, \ZZ_m$.
If $\Upsilon_\XX$ is an oriented cycle, then there exist $g \in \SSS_n$ and $m \in \{3, \dots, n\}$, such that $\XX=g\, \tilde\ZZ_m$.
\qed

\begin{figure}[ht]
\begin{center}
\begin{tabular}{ccc}
\parbox[c]{5.4cm}{\includegraphics[width=5.2cm]{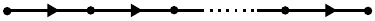}} 
& \quad &
\parbox[c]{5.4cm}{\includegraphics[width=5.2cm]{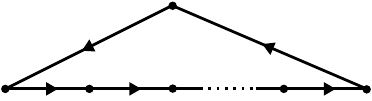}}
\end{tabular}

\bigskip
{\bf Figure 3.1.}
Oriented segment and oriented cycle.
\end{center}
\end{figure}

\bigskip\noindent
{\bf Proof of Proposition 2.3.}
Let $\XX$ be a full nonempty subset of $\SS$.
Write $\XX =\XX_1 \boxplus \cdots \boxplus \XX_\ell$, where $\XX_j$ is an indecomposable nonempty subset.
As observed above, each $\XX_j$ is also a full subset. 
Moreover, by Lemma 3.2, in order to show that $\varphi_\XX$ is faithful, it suffices to show that $\varphi_{\XX_j}$ is faithful for all $j \in \{1, \dots, \ell\}$.
So, we can assume that $\XX$ is a full and indecomposable nonempty subset of $\SS$.
By Lemma 3.3, there exist $g \in \SSS_n$ and $m \in \{2, \dots, n\}$ such that either $\XX = g\, \ZZ_m$, or $\XX = g\,\tilde \ZZ_m$ and $m \ge 3$.
Since $\varphi$ is equivariant under the actions of $\SSS_n$, upon conjugating by $g^{-1}$, we can assume that either $\XX = \ZZ_m$, or $\XX = \tilde \ZZ_m$.
Set $Z_m = \{x_1, \dots, x_m\} = \supp (\ZZ_m) = \supp (\tilde \ZZ_m)$, and identify $F_m$ with $F(Z_m)$.
Then $\varphi_{\ZZ_m} = \psi_m \circ f_m^{-1}$ and $\varphi_{\tilde \ZZ_m} = \tilde \psi_m \circ \tilde f_m^{-1}$, hence $\varphi_{\XX}$ is faithful by Theorem 3.1.
\qed


\section{Proof of Proposition 2.4}

The proof of Proposition 2.4 is based on some general results on Coxeter groups and Artin groups. 
Recall that the definitions of Coxeter diagram, Artin group and Coxeter group are given at the beginning of Step 2 in Section 2. 
Recall also that, if $Y$ is a subset of the set $S$ of vertices of $\Gamma$, then $\Gamma_Y$ denotes the full subdiagram spanned by $Y$, $A_Y$ denotes the subgroup of $A=A(\Gamma)$ generated by $Y$, and $W_Y$ denotes the subgroup of $W=W(\Gamma)$ generated by $Y$.

\bigskip\noindent
Let $\Gamma$ be a Coxeter diagram, let $S$ be its set of vertices, let $A$ be the Artin group of $\Gamma$, and let $W$ be its Coxeter group.
Since we have $s^2=1$ in $W$ for all $s \in S$, every element $g$ in $W$ can be represented by a word over $S$.
Such a word is called an \emph{expression} of $g$.
The minimal length of an expression of $g$ is called the \emph{length} of $g$ and is denoted by $\lg(g)$.
An expression of $g$ of length $\lg(g)$ is a \emph{reduced expression} of $g$.
Let $Y$ be a subset of $S$, and let $g \in W$.
We say that $g$ is \emph{$Y$-minimal} if it is of minimal length among the elements of the coset $W_Yg$.
The first ingredient in our proof of Proposition 2.4 is the following.

\bigskip\noindent
{\bf Proposition 4.1.}
(Bourbaki \cite[Chap. IV, Exercise 3]{Bourb1}).
{\it 
Let $Y \subset S$, and let $g \in W$. 
There exists a unique $Y$-minimal element lying in the coset $W_Y g$.
Moreover, the following conditions are equivalent.
\begin{itemize}[itemsep=2pt,parsep=2pt,topsep=2pt]
\item[(a)]
$g$ is $Y$-minimal,
\item[(b)]
$\lg (sg) > \lg (g)$ for all $s \in Y$, 
\item[(c)]
$\lg (hg) = \lg(h) + \lg(g)$ for all $h \in W_Y$.
\end{itemize}}

\bigskip\noindent
{\bf Remark.}
For $g \in W$ and $s \in S$, we always have either $\lg(sg) = \lg(g)+1$, or $\lg(sg) = \lg(g)-1$.
This is a standard fact on Coxeter groups that can be found for instance in \cite{Bourb1}.
So, the inequality $\lg(sg) > \lg(g)$ means $\lg(sg) = \lg(g)+1$ and the inequality $\lg(sg) \le \lg(g)$ means $\lg(sg) = \lg(g)-1$.

\bigskip\noindent
Let $u$ be a word over $S$.
\begin{itemize}[itemsep=2pt,parsep=2pt,topsep=2pt]
\item
Suppose that $u$ is written in the form $u_1 ss u_2$, where $u_1,u_2$ are words over $S$ and $s$ is an element of $S$. 
Then we say that $u' = u_1u_2$ is obtained from $u$ by an \emph{$M$-operation of type $\I$}.
\item
Suppose that $u$ is written in the form $u= u_1 \langle s,t \rangle^{m_{s,t}} u_2$, where $u_1,u_2$ are words over $S$ and $s,t$ are two elements of $S$ connected by an edge labelled by $m_{s,t}$.
Then we say that $u'= u_1 \langle t,s \rangle^{m_{s,t}} u_2$ is obtained from $u$ by an \emph{$M$-operation of type $\II$}.
\end{itemize}
We say that a word $u$ is \emph{$M$-reduced} if its length cannot be shortened by $M$-operations of type $\I,\II$. 
The second ingredient in our proof is the following.

\bigskip\noindent
{\bf Theorem 4.2}
(Tits \cite{Tits1}).
{\it Let $g \in W$.
\begin{itemize}[itemsep=2pt,parsep=2pt,topsep=2pt]
\item[(1)]
An expression $w$ of $g$ is a reduced expression if and only if $w$ is $M$-reduced.
\item[(2)]
Any two reduced expressions $w$ and $w'$ of $g$ are connected by a finite sequence of $M$-operations of type $\II$.
\end{itemize}}

\bigskip\noindent
Let $Y$ be a subset of $S$.
The third ingredient is a set-retraction $\rho_Y: A \to A_Y$ to the inclusion map $\iota_Y : A_Y  \to A$, constructed in \cite{GodPar1,ChaPar1}.
This is defined as follows.
Let $\alpha$ be an element of $A$.
\begin{itemize}[itemsep=2pt,parsep=2pt,topsep=2pt]
\item
Choose a word $\widehat \alpha = {s_1}^{\varepsilon_1} \cdots {s_\ell}^{\varepsilon_\ell}$ over $S^{\pm 1}$ which represents $\alpha$.
\item
Let $i \in \{0,1, \dots, \ell\}$.
Set $g_i = s_1s_2 \cdots s_i  \in W$, and write $g_i$ in the form $g_i = h_ik_i$, where $h_i \in W_Y$ and $k_i$ is $Y$-minimal.
\item 
Let $i \in \{1, \dots, \ell \}$.
If $\varepsilon_i =1$, set $z_i = k_{i-1}s_ik_{i-1}^{-1}$.   
If $\varepsilon_i =-1$, set $z_i = k_is_ik_i^{-1}$. 
\item
Let $i \in \{1, \dots, \ell \}$.
We set $T_i = z_i^{\varepsilon_i}$ if $z_i \in Y$.
Otherwise we set $T_i=1$.
\item
Set $\hat \rho_Y(\alpha) = T_1 T_2 \cdots T_\ell$.
\end{itemize}

\bigskip\noindent
{\bf Proposition 4.3}
(Godelle, Paris \cite{GodPar1}, Charney, Paris \cite{ChaPar1}). 
{\it Let $\alpha \in A$.
The element $\rho_Y(\alpha) \in A_Y$ represented by the word $\hat \rho_Y (\alpha)$ defined above does not depend on the choice of the representative $\widehat \alpha$ of $\alpha$.
Furthermore, the map $\rho_Y : A \to A_Y$  is a set-retraction to the inclusion map $\iota_Y : A_Y \hookrightarrow A$.}

\bigskip\noindent
We turn now to apply these three ingredients to our group $KB_n$ and prove Proposition 2.4.
Let $KW_n$ denote the quotient of $KB_n$ by the relations $\delta_{i,j}^2=1$, $1 \le i \neq j \le n$.
Note that $KW_n$ is the Coxeter group of the Coxeter diagram $\VGamma_n$.
For $\YY \subset \XX$, we denote by $KW_n(\YY)$ the subgroup of $KW_n$ generated by $\YY$.

\bigskip\noindent
{\bf Lemma 4.4.}
{\it Let $g \in KW_n$.
\begin{itemize}[itemsep=2pt,parsep=2pt,topsep=2pt]
\item[(1)]
An expression $w$ of $g$ is a reduced expression if and only if $w$ is $M$-reduced.
\item[(2)]
Any two reduced expressions $w$ and $w'$ of $g$ are connected by a finite sequence of $M$-operations of type $\II^{(2)}$ and $\II^{(3)}$.
\item[(3)]
Let $\YY$ be a subset of $\SS$, and let $w$ be a reduced expression of $g$. 
Then $g$ is $\YY$-minimal (in the sense given above) if and only if $w$ is $M_\YY$-reduced. 
\end{itemize}}

\bigskip\noindent
{\bf Proof.}
Part (1) and Part (2) are Theorem 4.2 applied to $KW_n$.
So, we only need to prove Part (3).

\bigskip\noindent
Suppose that $g$ is not $\YY$-minimal.
By Proposition 4.1, there exists $s \in \YY$ such that $\lg(sg)\le \lg(g)$, that is, $\lg(sg) = \lg(g)-1$.
Let $w'$ be a reduced expression of $sg$.
The word $sw'$ is an expression of $g$ and $\lg(sw') = \lg(w) = \lg(g)$, hence $sw'$ is a reduced expression of $g$.
By Theorem 4.2, $w$ and $sw'$ are connected by a finite sequence of $M$-operations of type $\II^{(2)}$ and $\II^{(3)}$.
On the other hand, $w'$ is obtained from $sw'$ by an $M$-operation of type $\III_\YY$.
So, $w'$ is obtained from $w$ by $M$-operations of type $\I$, $\II^{(2)}$, $\II^{(3)}$ and $\III_\YY$, and we have $\lg(w') < \lg(w)$, hence $w$ is not $M_\YY$-reduced.

\bigskip\noindent
Suppose that $w$ is not $M_\YY$-reduced. 
Let $w'$ be an $M_\YY$-reduction of $w$, and let $g'$ be the element of $KW_n$ represented by $w'$.
Since $w'$ is an $M_\YY$-reduction of $w$, the element $g'$ lies in the coset  $KW_n(\YY)\,g$.
Moreover, $\lg(g') = \lg(w') < \lg(w) = \lg(g)$, hence $g$ is not $\YY$-minimal.
\qed

\bigskip\noindent
{\bf Proof of Proposition 2.4.}
Let $\YY$ be a subset of $\SS$.
Consider the retraction $\rho_\YY : KB_n \to KB_n(\YY)$ constructed in Proposition 4.3.
We shall prove that, if $u$ is a word over $\SS^{\pm 1}$, then $\overline{\pi_\YY(u)} = \rho_\YY (\bar u)$.
This will prove Proposition 2.4.
Indeed, if $\bar u = \bar v$, then $\overline{ \pi_\YY (u)} = \rho_\YY (\bar u) = \rho_\YY (\bar v) = \overline{\pi_\YY (v)}$.
Moreover, since $\rho_\YY : KB_n \to KB_n(\YY)$ is a retraction to the inclusion map $KB_n(\YY) \hookrightarrow KB_n$,  we have $\rho_\YY (\bar u) = \bar u$ if and only if $\bar u \in KB_n(\YY)$, hence $\overline{ \pi_\YY (u)}  = \bar u$ if and only if $\bar u \in KB_n(\YY)$.

\bigskip\noindent
Let $u = s_1^{\varepsilon_1} \cdots s_\ell^{\varepsilon_\ell}$ be a word over $\SS^{\pm 1}$.
Let $\alpha$ be the element of $KB_n$ represented by $u$.
\begin{itemize}[itemsep=2pt,parsep=2pt,topsep=2pt]
\item
For $i \in \{0,1, \dots, \ell\}$, we set $u_i^+=s_1 \cdots s_i$, and we denote by $g_i$ the element of $KW_n$ represented by $u_i^+$.
\item
Let $i \in \{0,1, \dots, \ell\}$.
We write $g_i = h_i k_i$, where $h_i \in KW_n(\YY)$, and $k_i$ is $\YY$-minimal.
Let $v_i^+$  be an $M_\YY$-reduction of $u_i^+$.
Then, by Lemma 4.4, $v_i^+$ is a reduced expression of $k_i$. 
\item
Let $i \in \{1, \dots, \ell\}$.
If $\varepsilon_i =1$, we set$z_i = k_{i-1} s_i k_{i-1}^{-1}$ and $w_i^+=v_{i-1}^+ \cdot s_i \cdot \op(v_{i-1}^+)$.
If $\varepsilon_i =-1$, we set $z_i = k_i s_i k_i^{-1}$ and $w_i^+=v_i^+ \cdot s_i \cdot \op(v_i^+)$.
Note that $w_i^+$ is an expression of$z_i$.
\item
Let $i \in \{1, \dots, \ell\}$.
Let $r_i$ be an $M$-reduction of $w_i^+$.
By Lemma 4.4, $r_i$ is a reduced expression of $z_i$.
Note that we have $z_i \in \YY$ if and only if $r_i$ is of length $1$ and $r_i \in \YY$.
\item
Let $i \in \{1, \dots, \ell \}$.
If $r_i$ is of length $1$ and $r_i \in \YY$, we set $T_i = r_i^{\varepsilon_i}$.
Otherwise we set $T_i=1$.
\item
By construction, we have $\hat \rho_\YY(\alpha) = \pi_\YY(u) = T_1 T_2 \cdots T_\ell$.
\end{itemize}
\qed



\bigskip\bigskip\noindent
{\bf Paolo Bellingeri}, 

\smallskip\noindent
LMNO UMR6139, CNRS, Université de Caen, F-14000 Caen, France.

\smallskip\noindent
E-mail: {\tt paolo.bellingeri@unicaen.fr}

\bigskip\noindent
{\bf Bruno A. Cisneros de la Cruz},

\smallskip\noindent 
Instituto de Matemáticas de la UNAM - Oaxaca, Oaxaca de Juárez, Oax. 68000, Mexico.

\smallskip\noindent
E-mail: {\tt brunoc@matem.unam.mx}

\bigskip\noindent
{\bf Luis Paris},

\smallskip\noindent 
IMB UMR5584, CNRS, Univ. Bourgogne Franche-Comté, F-21000 Dijon, France.

\smallskip\noindent
E-mail: {\tt lparis@u-bourgogne.fr}

\end{document}